\documentclass[conference]{IEEEtran}
\bstctlcite{IEEEexample:BSTcontrol}

\usepackage{lineno,hyperref}
\usepackage{amsmath,amssymb,fixmath}
\usepackage{graphicx}

\begin{document}
\title{Quantum Calculus-based Volterra LMS for Nonlinear Channel Estimation}


 \author{\IEEEauthorblockN{Muhammad~Usman\IEEEauthorrefmark{1}, Muhammad~Sohail~Ibrahim\IEEEauthorrefmark{1}, Jawwad~Ahmad\IEEEauthorrefmark{2}, Syed~Saiq~Hussain\IEEEauthorrefmark{4}, and Muhammad~Moinuddin\IEEEauthorrefmark{3}}
 	\IEEEauthorblockA{\IEEEauthorrefmark{1}Iqra University, Karachi, Pakistan.\\ Email: {musman, msohail}@iqra.edu.pk\\ 		
 		\IEEEauthorrefmark{2}Usman Institute of Technology, Karachi-75300, Pakistan.\\
 		Email: {jawwad}@uit.edu\\
    	\IEEEauthorrefmark{4}College of Engineering, PAF-KIET, Karachi, Pakistan.\\
 		Email: {saiqhussain}@gmail.com\\
        \IEEEauthorrefmark{3}CEIES, King Abdulaziz University, Saudi Arabia.\\
 		Email: {mmsansari}@kau.edu.sa\\ 		
 		}}
\maketitle
\begin{abstract}
A novel adaptive filtering method called $q$-Volterra least mean square ($q$-VLMS) is presented in this paper.  The $q$-VLMS is a nonlinear extension of conventional LMS and it is based on Jackson's derivative also known as $q$-calculus.  In Volterra LMS, due to large variance of input signal the convergence speed is very low. With proper manipulation we successfully improved the convergence performance of the Volterra LMS.  The proposed algorithm is analyzed for the step-size bounds and results of analysis are verified through computer simulations for nonlinear channel estimation problem.
\end{abstract}


\IEEEpeerreviewmaketitle
\section{Introduction}\label{Sec:Intro}
With the advent of high flexible bandwidth and mobility, the modern wireless communication systems are able to provide high throughput values. These inventions however have to deal with numerous challenges. One of the utmost requirement of an efficient communication system is to be able to effectively perform channel estimation process. Linear models provide reasonable estimate of the channel alongwith guarantee reduced Bit Error Rate (BER). Therefore these models are considered to play significant role in the improvement of system's performance \cite{FLMF}. 

The concept of fractional calculus has been widely incorporated in various research areas \cite{FBPTT, FRBF, Sohail, Usman}. In \cite{FCLMS}, fractional order calculus was utilized a to propose a least mean square (LMS) algorithm. Complex linear systems were successfully identified using the algorithm. It is reported that the proposed FCLMS algorithm provides better convergence rate than the conventional complex least mean square (CLMS) algorithm.

However, some situations demand the utilization of nonlinear solutions \cite{VLMS_MOC,Bilinear_VLMS}. In \cite{RVP, RVSS}, fractional least mean square was utilized to propose an adaptive framework with variable power. The proposed method was applied on channel equalization and plant identification problems and it is shown that the algorithm adopts the fractional power. This dynamically adds the levels of derivatives and in turn we can achieve a better convergence rate. Another important parameter is error of the steady state, which can also be kept at minimum by utilizing the fractional power. It was reported in  \cite{VPFLMS}, that the instantaneous error energy aids the adaptation of the fractional power of fractional gradient based least mean square to obtain faster convergence and lower errors. However, the fractional variants of LMS are inherently prone to some divergence issues, as their designs are not fully mathematically plausible and extensive study of convergence and stability analysis is missing.  In addition to this, numerous studies has already been conducted, highlighting the flaws in the design and implementation issues of FLMS and other fractional variants on different problems of engineering and scientific interest. For instance in \cite{Comments1}, authors discussed that the fractional order calculus variants of least mean square are unstable and the design is invalid for negative weights and complex signals. In \cite{comments2}, author highlighted the major flaws in the design and implementation of normalized FLMS and its extension to complex signal processing. Similarly, in \cite{comments3}, another variant of FLMS is criticised for its pseudo improved performance gain, when the original design of FLMS is modified by introducing the momentum term and it is implemented for the identification of power signal parameters. In \cite{comments3}, it is found that the performance of FLMS and its momentum based variant is not higher than the conventional LMS and in some cases it is inferior to the standard LMS algorithm. Due to these reasons many authors are quite skeptical in using FLMS algorithm for real world problems.

To design a novel state-of-the art variant of LMS using the notion of novel gradient descent approach, many researchers look forward towards Jackson derivative-based LMS which is also known as q-calculus based LMS or simply q-LMS \cite{eqLMS}. 
In \cite{eqLMS}, an enhanced q-calculus based LMS was proposed that incorporated time varying q-parameter that utilizes parameter-less concept of error-correlation energy to produce higher stability, convergence, and lower steady state errors. Another variant of q-calculus based LMS was proposed in \cite{QLMF}, where q-calculus based least mean fourth was implemented to cater the non-Gaussian noise based channel estimation problem.

A simple solution to deal with the nonlinearity is the Volterra LMS, which is a nonlinear extension of LMS using Volterra series expansion. In this study we tend to use quantum calculus based volterra LMS and observe its efficacy on non-linear channel estimation problem. The remainder of the paper provides following details: An introduction to the volterra filter is presented in section \ref{volterra}, q-VLMS is analyzed in Section \ref{Sec:Analysis}. The results for the experiments are discussed in Section \ref{Sec:Sim} followed by conclusion in \ref{Sec:Con}.

\section{Volterra Filter} \label{volterra}

A second order expansion of Volterra filter can be defined as \cite{VLMS} :
\begin{multline} \label{Volterra_system}
y(r) = h_{0} + \Sigma_{d=0}^{M-1} A(d) x(r-d) \\ + \Sigma_{d=0}^{M-1}\Sigma_{e=0}^{M-1} B(d,e) x(r-d) x(r-e)
\end{multline}
where ${A(d)}$ represnts the weight for linear filter while  ${B(d,e)}$ represents quadratic filter weights. $M$ denotes the filter length. The relationship can also be written as: 
\begin{equation}
y(r) = \mathbf{w}^{\intercal}(r) \mathbf{u}(r),
\end{equation}
where, $\mathbf{u}(r) = [u_1,u_2, \dots, u_N] = [x(r), x(r-1), \dots, x(r-M+1), x(r)x(r), x(r)x(r-1), \dots, x(r-M+1)x(r-M+1)]$,  and $\mathbf{w}(r)= [w_1,w_2, \dots, w_N]=[a(1), a(2), \dots, A(M), B(1,1), B(1,2), \dots B(M,M)]$.

Using the conventional gradient descent method:
\begin{equation} \label{VLMS_grad}
\boldsymbol{\Delta w}(r) =  - \frac{\mu}{2} \nabla_{\mathbf{w}}C(\mathbf{w}),
\end{equation}

where $C(\mathbf{w})=E[e^2(r)]$ represents the average loss function, error $e(r)=o(r)-y(r)$ is calculated by taking difference of output $y(r)$ and desired signal $o(r)$ at $r$th iteration.
The weights update rule for the $l$th weight of the Volterra LMS filter is derived as
\begin{equation*} \label{VLMS_WU}
w_{l}(r+1) = w_{l}(r) - \mu * e(r) * u(r),
\end{equation*}
where $\mu$ is the step-size of the gradient descent algorithm.  The input correlation matrix has a huge influence on the performance of the Volterra LMS and it is highly dependent on the eigenspread of the input signal.  To resolve this underlying issue q-derivative is used to design a non-linear variant of VLMS. Instead of conventional gradient descent approach, we propose to use the Jackson derivative method \cite{eqLMS}, which takes steps towards the optima at higher speed i.e. (for $q>1$) the algorithm optimize the adaptive parameters rapidly in the search direction. An example of which can be found in \cite{eqLMS} where instead of utilizing the conventional tangent of the cost function, secant is used. Consequently by substituting q-gradient based novel stochastic gradient descent algorithm in place of the conventional gradient descent algorithm in \eqref{VLMS_grad}, we get:
\begin{eqnarray}\label{q_gradient}
\boldsymbol{w}(r+1) = \boldsymbol{w}(r) - \frac{\mu}{2} \nabla_{q,w}C(w).
\end{eqnarray}
The q-derivative can be calculated as $d_{q}(f(r)) = f(qn)-f(r)$. Eventually it is turns out to be:
\begin{equation}
D_{q}(f(r)) = \frac{d_{q}(f(r))}{d_{q}(x)} = \frac{f(qn)-f(r)}{(q-1)i}.
\end{equation}
if $q$ is chosen as $1$ the derivative is considered as a conventional derivative. Applying the $q$-gradient the cost function $C(\mathbf{w})=E[e^2(r)]$  becomes:
\begin{eqnarray}\label{q_gradient2}
\nabla_{q,w}C(w) = -2E[\mathbf{G} \boldsymbol{u}(r) e(r)]
\end{eqnarray}
where 
${\rm diag}(\mathbf{G}) = [(\frac{q_{1}+1}{2}), (\frac{q_{2}+1}{2}),.....(\frac{q_{M}+1}{2})]^{\intercal}$. 

By dropping the expectation in \eqref{q_gradient2} results in $\nabla_{q,w}C(w) \approx -2G \boldsymbol{u}(r) e(r).$ which upon substitution in \eqref{VLMS_grad} deduce to q-VLMS algorithm:
\begin{eqnarray}\label{qLMS_final}
\boldsymbol{w}(r+1) = \boldsymbol{w}(r) +\mu G \boldsymbol{u}(r)e(r).
\end{eqnarray}

\section{Analysis}\label{Sec:Analysis}
\subsection{Optimal Solution for the $q$-VLMS}
Consider $\mathbf{u}$ as Gaussian white noise having a unit variance and a mean value of zero.  Based on \cite{R_diagonalization,VLMS_Analysis}, the input vector for the adaptive plant is $\mathbf{S}^{-1} \mathbf{u}(r) = \mathbf{S}^{-1}[u_1,u_2, \dots, u_M]^{\intercal}$, where $\mathbf{S}^{-1}$ as shown in equation \ref{seqn},  serves the purpose to scale the square minus one terms by $\frac{1}{2}$.
\begin{eqnarray}\label{seqn}
diag(\mathbf{S}) = \{1,1,1,\sqrt[]{2},1,1,\sqrt[]{2},1,\sqrt[]{2}\} .
\end{eqnarray}
The purpose of introducing the $\mathbf{S}^{-1}$ is to make input auto-correlation matrix become an identity matrix i.e, $\mathbf{S}^{-1} \boldsymbol{R}\mathbf{S}^{-1}= \boldsymbol{I}$. The output of second-order VLMS is now defined as: $y(r) = \mathbf{w}^{\intercal}(r) \mathbf{S}^{-1} \mathbf{u}(r)$.

For optimal solution replacing $\nabla_{q,w}C(w)$ in \eqref{q_gradient2} with zero.
\begin{eqnarray}\label{q_gradient3}
\mathbf{G} E[\mathbf{S}^{-1} \boldsymbol{u}(r) e(r)] =0,
\end{eqnarray}
\begin{eqnarray}\label{q_gradient4}
E[ \mathbf{S}^{-1} \boldsymbol{u}(r) \boldsymbol{d}(r)] -E[ \mathbf{S}^{-1} \boldsymbol{u}^{\intercal}(r) \boldsymbol{u}(r)\mathbf{S}^{-1}] E[\boldsymbol{w}(r)] =0,
\end{eqnarray}
\begin{eqnarray}\label{q_gradient5}
\mathbf{S}^{-1} \boldsymbol{r}_{ud} - \mathbf{S}^{-1}\boldsymbol{R}\mathbf{S}^{-1} \boldsymbol{w}_{opt} =0,
\end{eqnarray}
where $\boldsymbol{R}$ and $\boldsymbol{r}_{ud}$ represents the input auto-correlation matrix and the desired output and input vector of cross-correlation respectively.

To obtain the optimal weight vector we use $ \boldsymbol{w}_{opt} = \mathbf{S} \boldsymbol{R}^{-1} \mathbf{S}  \mathbf{S}^{-1} \boldsymbol{r}_{ud} = \mathbf{S} \mathbf{w}^{*}$, where
\begin{eqnarray}\label{q_gradient6}
\boldsymbol{w}^{*} =\boldsymbol{R}^{-1} \boldsymbol{r}_{ud} ,
\end{eqnarray}
which is same as optimal Wiener solution of the LMS, and the minimum least square error at $\boldsymbol{w}_{opt}$ is similar to Wiener power if linear LMS is considered.
\begin{eqnarray}\label{q_gradient7}
\boldsymbol{\xi}_{min} = E[\eta^{2}(r)],
\end{eqnarray}
where $\eta$ is the Gaussian noise.

\subsection{Convergence Analysis}
In this section we analyse the proposed  $q$-VLMS  for both mean and mean square performances with the following commonly used assumptions \cite{VLMS_Analysis}: (1) The input and noise signals are random signals have normal Gaussian distribution i.e unit variance and zero mean, (2) the input sequence vector $\mathbf{u}$ is i.i.d i.e independent and identically distributed.

The weight error vector is defined as $\Delta_{w}(r)=\boldsymbol{w}_{opt}-\boldsymbol{w}^{*}(r)$, $
e(r)=\Delta_{w}^{\intercal}(r)\boldsymbol{u}(r) + \eta(r)$. After substituting $\boldsymbol{w}(r)=\boldsymbol{w}^{*}(r)$, $\boldsymbol{u}(r)=\boldsymbol{S}^{-1} \boldsymbol{u}(r)$, and $e(r)$ in \eqref{qLMS_final}, we get
\begin{equation}\label{Mean_analysis1}
\Delta_{w}(r+1) = \Delta_{w}(r) +\mu \boldsymbol{G} \boldsymbol{S}^{-1} \boldsymbol{u}(r)( \boldsymbol{u}^{\intercal}(r)\boldsymbol{S}^{-1} \Delta_{w}(r) + \eta(r)).
\end{equation}
which upon simplification results in:
\begin{eqnarray}\label{Mean_analysis3}
E\left[\Delta_{w}(r)\right] = \left(\mathbf{I} - \mu \mathbf{A}\right)^{i} E\left[\Delta_{w}(0)\right].
\end{eqnarray}
where $\mathbf{A} = \boldsymbol{G} E\left[\boldsymbol{S}^{-1} \boldsymbol{u}(r)( \boldsymbol{u}^{\intercal}(r)\boldsymbol{S}^{-1}\right]$ and $E[]$ is the expectation operator.

For convergence 
\begin{equation}
0 < \mu < \frac{1}{max\{(q_1 +1)\lambda_{1},\dots,(q_M +1)\lambda_{M} \}}
\end{equation}
In case when all $q_r$'s are equal to $q$, the limits will be $0 < \mu < \frac{1}{(q+1)\lambda_{max}}$, where $\lambda_{max}$ is the maximum eigenvalue of $S^{-1}\boldsymbol{R}S^{-1}$.

\subsection{Computational Complexity}
The proposed $q$-VLMS is just $K$ multiplication expensive than the conventional VLMS.  It takes $3K+1$ multiplications because in q-gradient we also need to calculate the multiplication of  diagonal matrix $\boldsymbol{G}$. In addition to the above mentioned computation we also need to perform $2K$ additions, one for the calculation of error signal and one for the weight update.  Here $K$ symbol represents the total number of adaptive parameters.

\section{Experiments}\label{Sec:Sim}
To evaluate the performance of proposed nonlinear variant of LMS also known as q-Volttera LMS, consider a nonlinear channel as defined in \eqref{Volterra_system} of order (N=3), $d(r) = \mathbf{h}^{\intercal}(r) \mathbf{u}(r) + \eta(r)$, where $u(r)$ is the input signal which is a vector of length equal to the length of weight vector  and $d(r)$ is referred to the output of the system. Noise signal is represented by $\eta(r)$, which in system identification term also known as the  disturbance (noise) taken to be white noise generated using a randomly distributed Gaussian source.  For the experiment, input signal is chosen to be a Gaussian distributed source having variance of $1$ and zero mean.  

For the performance comparison and evaluation on standard measures, we opt to compare the estimated weights with the  actual ones  using normalized weight deviation (NWD) method. In particular,  it is defined as $\mbox{NWD}=\frac{\left\Vert \bf h-\bf w \right\Vert^{2}}{\left\Vert\bf h\right\Vert^{2}}$.  For fair evaluation $1000$ independent trails of the simulations are performed and average result is reported.  For each simulation round the coefficient values of the desired channel and the initial weights of the adaptive filter were randomly selected.

To validate the analysis results and compare the performance of the proposed algorithm with contemporary methods. We design two evaluation protocols:

\begin{enumerate}
	\item {\bf Evaluation protocol 1}: Validation of analytical results on three values of q=$\{1, 5, 10\}$, at signal-to-noise ratio (SNR) of 20dB, and $\mu=0.25\times \lambda_{\max}^{-1}$.
	\item {\bf Evaluation protocol 2}: Effect of $q$-parameters amd performance comparison with conventional VLMS on learning rate $\mu=1\times10^{-3}$, and three noise levels SNR=\{10, 20 ,30\} dB.
\end{enumerate}

 \begin{figure}[htb]
 	\centering
 	\centerline{\includegraphics[width=0.45\textwidth]{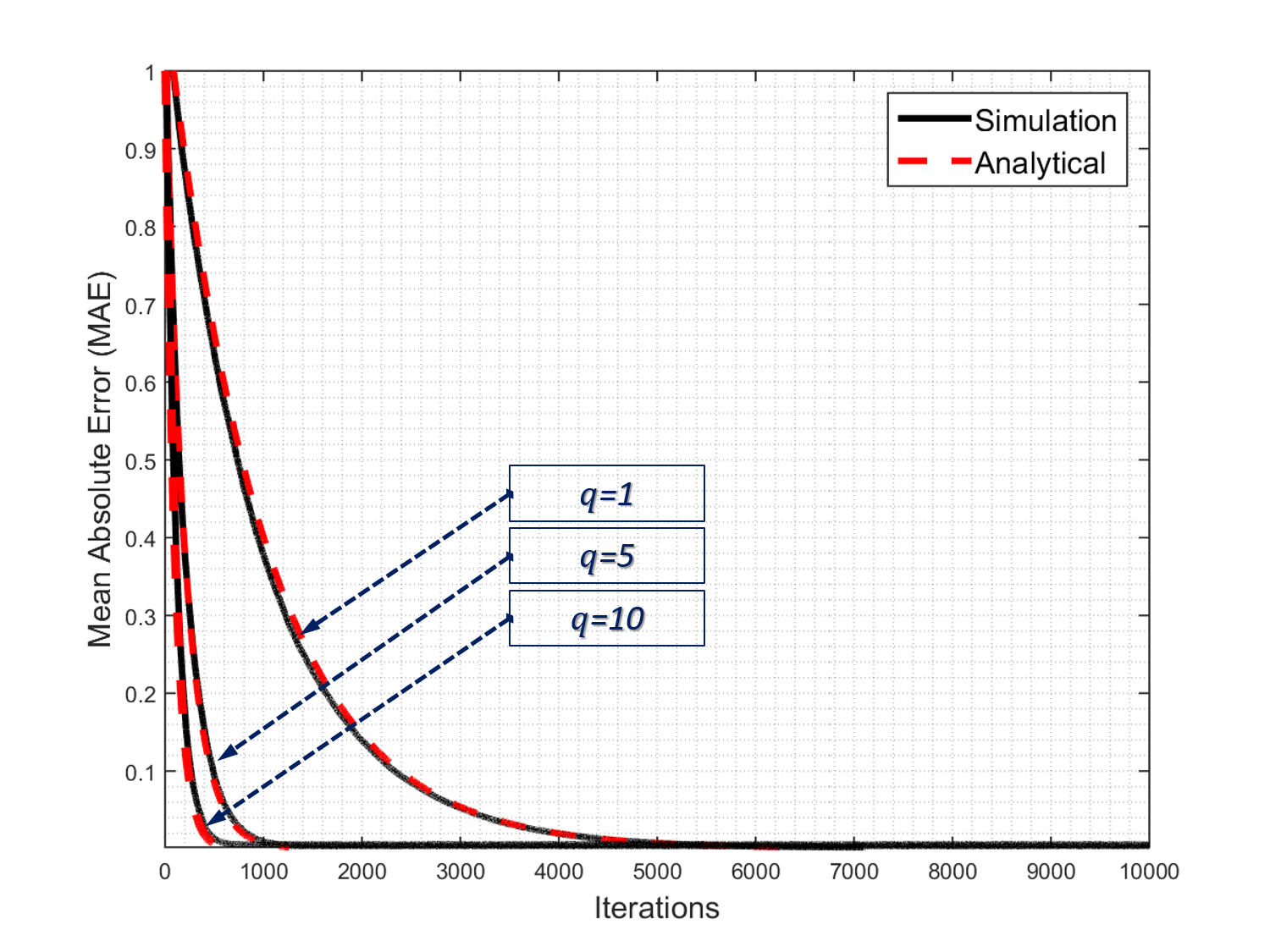} }
 	\caption{Validation of analytical results.}
 	\label{fig:resultsCon}
 \end{figure}
 
Fig. \ref{fig:resultsCon} show the comparison of analytical and simulation results for mean absolute error at three different values of $q$.  We choose to verify three different values of q i.e. $q=\{1, 5, 10\}$. Interestingly, the analytical results are well matched with simulations and in order to quantitatively measure the similarity between two results we measure the correlation coefficient of MAE between simulation and analysed values, and the results are well correlated with an average correlation coefficient of $0.9995$.  

Using evaluation protocol 2, the sensitivity of $q$ parameter is analyzed for three different noise levels. In particular, we implemented system identification problem with signal-to-noise ratios of 10, 20 and 30 dB. The proposed algorithm is also compared with the VLMS, we repeat the same simulations with $G=\mathbf{S}\boldsymbol{R}^{-1}\mathbf{S}$.
 
 \begin{figure}[htb]
 	\centering
 	\centerline{\includegraphics[width=0.45\textwidth]{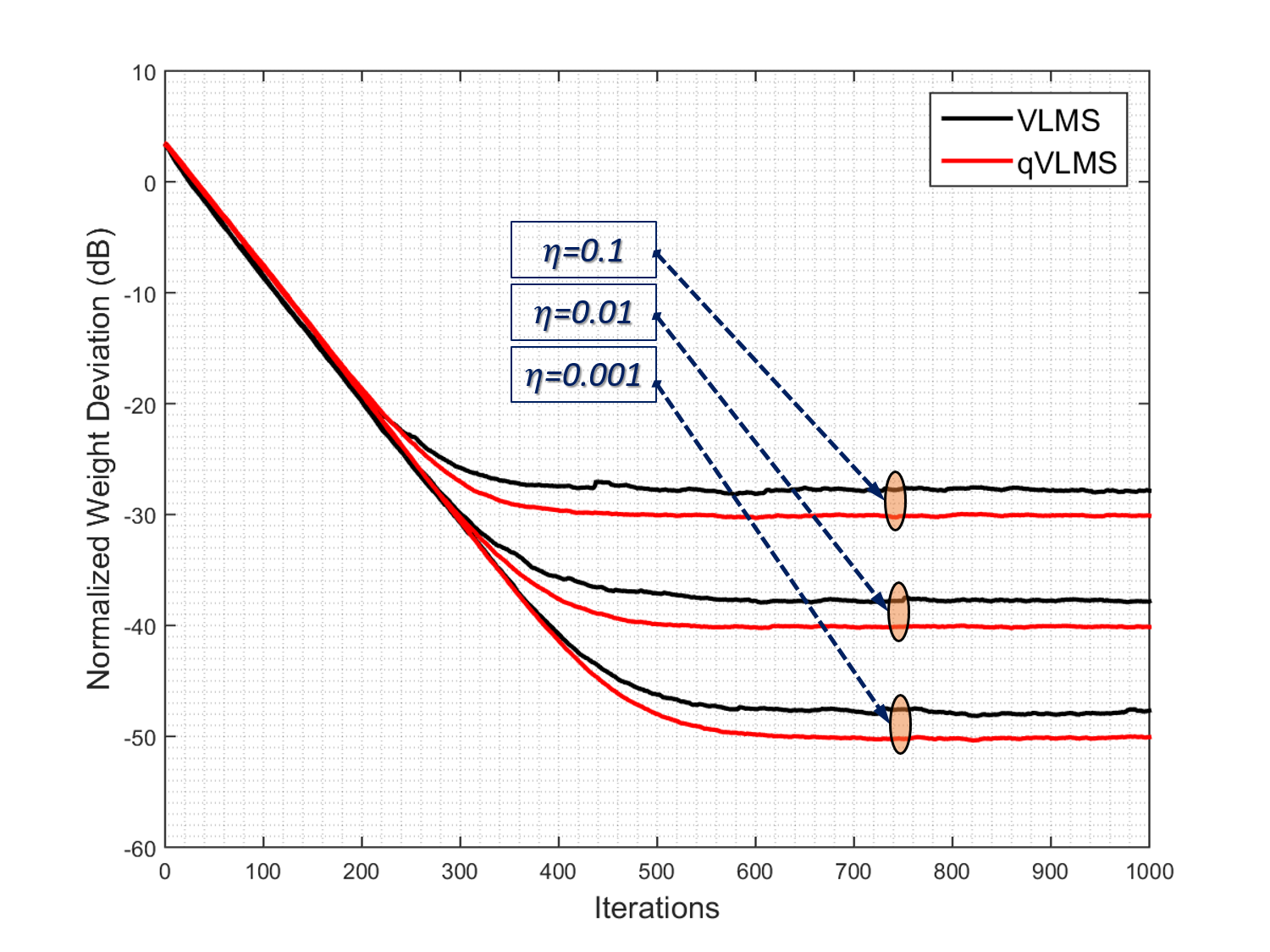} }
 	\caption{Effect of $q$-parameters.}
 	\label{fig:convergence_scenario}
 \end{figure}

The purposed $q$-VLMS achives an excellent performance gain for all given scenarios. It is evident from Fig. \ref{fig:convergence_scenario}, that the the NWD values of the proposed method are superior in for all SNRs. It outperformed the conventional VLMS by an average $2.31$ dB.

\section{Conclusion}\label{Sec:Con}
In this study, we presented a $q$-calculus-based LMS algorithm names as $q$-VLMS.  The proposed method provide additional control over the convergence and steady-state performances through $q$ adjustment parameters.  The proposed $q$-VLMS was compared with the conventional VLMS algorithm for the estimation of non linear channel. The performance of the algorithm was tested for different values of SNRs and convergence and steady state values were observed. The algorithm outperformed the conventional VLMS by achieving better results in all test scenarios and evaluation parameters.



\end{document}